\documentclass{amsart}
\usepackage{graphics,verbatim,color,amsthm}

\usepackage[all]{xy} 
\theoremstyle{plain}
\newtheorem{theorem}{Theorem}

\theoremstyle{definition}

\def\R{\mathbb{R}}

\def\C{\mathbb{C}}

\def\bP{\mathbb{P}}

\def\V{\mathbb{V}}

 \newcommand{\dd}[2]
{
{{\partial #1}   \over {\partial #2}}
}

\begin{document}

\author{Richard Montgomery}
\address{Mathematics Department\\ University of California, Santa Cruz\\
Santa Cruz CA 95064}
\email{rmont@ucsc.edu}

\date{September 15, 2016}

\title{Constructing the  Hyperbolic Plane as the reduction of  a three-body problem}

\begin{abstract} We construct  the hyperbolic plane with its geodesic flow as
the scale plus symmetry reduction of a   three-body problem in the Euclidean plane.
The potential is $-I/\Delta^2$ where $I$ is the triangle's moment of inertia and $\Delta$ its area.
The  reduction method uses the Jacobi-Maupertuis metric, following   \cite{hyperbolic_pants}. 
\end{abstract} 

\maketitle

%Inspired in part by the inverse problem solved by Fujiwara et al
%and by a conversation with Persi Diaconis, 

\section{The problem.} Three  point particles  move in the Euclidean plane $\R^2 = \C$  according to Newton's equations
\begin{equation}
m_a \ddot q_a = - \nabla_a V, a=1,2,3, 
\end{equation}
with potential  
\begin{equation}
V(q) = -\frac{ \gamma \text{ (moment of inertia)} }{(\text{area})^2} = -\frac{\gamma I (q)}{ \Delta (q) ^2}.
\label{potential}
\end{equation}
The $m_a$ are the   masses.  The $q_a \in \R^2 = \C$ are the instantaneous positions of these point masses.  
The   $\nabla_a = \dd{}{q_a}$  are the    gradients with respect to   $q_a$.  We write   $q = (q_1, q_2, q_3) \in \C^3$ for the
three positions put together into one vector. 
The denominator of our potential is
$$\Delta (q) = \text{ signed area of triangle of triangle with vertices } q_1, q_2, q_3,$$
% is the  signed area of the triangle whose vertices are   $q_1, q_2, q_3$. 
 while its   numerator  is  
$$ I(q)=\text{ moment of inertia} = \frac{\Sigma_{a <b} m_a m_b r_{ab} ^2 } {\Sigma_a m_a} , \text{ where }  r_{ab} = |q_a - q_b|, $$
 the moment of inertia of this triangle with  respect to its center of mass. 
The constant  $\gamma > 0$ is a physical constant needed to make   the units of the potential  that of energy, so that $\gamma$
has  units of $(\text{length})^4/\text{(time)}^2$.  

Our ODEs form a Galilean-invariant Hamiltonian system, and as such have the usual conserved quantities  
$$H = K(v) + V(q),  P = \Sigma m_a v_a,   J = \Sigma m_a q_a \wedge v_a , $$
of energy, linear momentum and angular momentum.  In the  expression for   energy the term   $K$ is the kinetic energy
$$K(v) =  \frac{ 1}{2}  \Sigma m_a |v_a|^2 = \frac{1}{2} \langle v, v \rangle, $$
a function of   the velocities $v = \dot q \in \C^3$.  The inner product occuring here is the
``mass inner product'' on $\C^3 = (\R^2)^3$,  $$\langle q, v \rangle = m_1 q_1 \cdot v_1  + m_2 q_2 \cdot v_2  + m_3 q_3 \cdot v_3 .$$
%This Hamiltonian  system  is Galilean-invariant, so has a conserved linear momentum $\Sigma m_a v_a$  and angular momentum $J = \Sigma m_a q_a \wedge v_a$.  
where the dot product $q_i \cdot v_i = Re (q_i \bar v_i)$ is the usual dot product in $\R^2 \cong \C$.  Standard physics tricks tell us that it is no loss of generality to restrict ourselves to the   center-of-mass subspace:
$$\V_{cm} = \{q \in \C^3: \Sigma m_a q_a = 0 \} .$$ 

On the center-of-mass subspace  
 Lagrange showed that 
\begin{equation}
\label{Lag}I = \langle q, q \rangle : =\Sigma m_a |q_a|^2  
\end{equation} 
A standard computation, yields the Lagrange-Jacobi (or virial) identity
$$\ddot I = 4H$$
which is valid for any potential $V$ which is homogeneous of degree $-2$. 
It follows that on the energy level $H=0$, the phase space function $\dot I  =2 \langle q, v \rangle $ is an
additional concerved quantity.   If we do not want the size $I$  of our solution to change then, we   work on the
invariant submanifold of phase space for which   $H=0, \dot I = 0$.  

Our problem is to solve our Newton's equations on the  submanifold of
phase space for  which $H =0, P=0, \dot I = 0$ and $J=0$.

\section{Reduction and solution.} 
The  group $G$ of rigid motions  of the plane is a subgroup of the Galilean group and so   maps solutions to solutions.
As a consequence, the dynamics pushes down to  the quotient space of $\C^3$ by $G$.
This quotient space is  called ``shape space'' and is  homeomorphic to $\R^3$. Points of shape space are oriented congruence classes of triangles.  Shape space
is  endowed with standard ``Hopf-Jacobi' coordinates $w_1, w_2 ,w_3$ which correspond to  quadratic $G$-invariant functions on $\C^3$.
These coordinates have the property that  $$w_3 = \mu  \Delta; \mu ^2 = \frac{m_1 m_2 m_3}{m_1 + m_2 + m_3} $$ 
and $w_1 ^2 + w_2 ^3 + w_3 ^2 = \frac{I^2}{4}$.   See \cite{me_monthly} or
\cite{remarkable}.    Note that the equator $w_3 = 0$ corresponds to the degenerate triangles in which all three masses lie along a single line.  We call $\pi: \C^3 \to \R^3$; $q \mapsto \pi(q) = w$
the ``shape space projection''.   ( The  formal process of pushing down the dyanmics is   ``symplectic reduction''
and requires fixing the value of the angular momentum $J$.  We have already fixed the linear momentum $P$ to be zero and will 
be fixing the angular momentum to be zero.)

Within shape space is the {\it shape sphere} $S^2$, which is the sphere $|w|^2 = 1$ (or any sphere $I = I_0$ as long as $I_0$ is a positive   constant).
Points of this sphere are identified with oriented similarity classes of triangles,
since scaling a triangle $q$  by $\lambda>0 $ multiplies $I$ by $\lambda^2$.    
Alternatively, if we delete the origin from $\V_{cm} \cong \C^2$ and then quotient by complex scaling
$q \mapsto \lambda q$, $\lambda \ne 0$ in $\C$, we arrive at $\C \bP^1 \cong S^2$.
Now reflection about any line
in the plane of the triangles has the effect $(w_1, w_2, w_3) \mapsto (w_1, w_2, - w_3)$
so that the upper hemisphere $w_3 \ge 0$ of the shape sphere $|w|^2 = 1$,  realizes the space of similarity classes of triangles, with its boundary $w_3 = 0$ representing the degenerate collinear triangles.

The open upper hemisphere, $w_3 > 0$, $|w|^2 = \frac{I_0 ^2}{4} $, endowed with the     metric $\frac{dw_1 ^2 + dw_2 ^2 + dw_3 ^2}{w_3 ^2}$   
is a fairly  well-known realization of the hyperbolic plane, sometimes referred to as the ``Jemisphere model''.  See pp. 69-71 of \cite{All_Models}. 
%  Thus the Jemisphere metric is a metric on the space of (nondegenerate!) similarity classes of triangles. 

\begin{theorem}  The shape space projection of those  solutions to our Newton's equations  for which $H= \dot I = J = P = 0$ 
are  geodesics  for the Jemisphere model of the   hyperbolic plane.   Relative to the standard Hopf-Jacobi  induced Cartesian  coordinates $w_1, w_2, w_3$ on
shape space described above, 
these projected solution curves are obtained by intersecting a hemisphere 
$w_1 ^2 + w_2 ^2 + w_3 ^2 = (1/4) I_0 ^2$ ,   $w_3  > 0$ or $w_3 < 0$ 
with  a vertical plane $A w_1 + Bw_2 = const. $. (See figure \ref{Klein_to_Jemi}.)
\end{theorem}

{\bf Remark 1.}  The reflection $(w_1, w_2, w_3) \mapsto (w_1, w_2, - w_3)$ maps the ``upper Jemisphere model'' $w_3 > 0$
canonically isometrically to the lower Jemisphere model $w_3 < 0$, both having the same form above for the metric.

{\bf Remark 2.} As a consequence of the expression of the projected solutions $w(t) = \pi (q (t))$, we see that when such  a solution
is extended over its
maximum time range $(a,b)$ of existence, it begins at one collinear configuration and ends at another. (In particular  $lim_{t \to a, b} w_3(t) = 0$.)  
One can  compute the angle between the   lines containing these two extremal collinear configurations by  using the area formula,
as per the periodicity proof for the figure eight solution given in \cite{remarkable}.

\section{Proof}  

%The proof is essentially a computation. 
%We begin with some

% {\it  Computational Fundamentals. }
%We first discuss the `scale reduction' which is only available for potentials which are
%homogeneous of degree $-2$. 

{\bf Step 1. Jacobi-Maupertuis.}   The  Jacobi-Maupertuis principle, applied for energy   $H=0$,    asserts that
the zero-energy solutions to our Newton's equations are, up to reparameterization,   geodesics for the metric  
$$ds^2_{JM} = 2 U(q) |dq|^2_E , \text{ where } U = -V$$
on $\C^3$.   The center-of-mass zero subspace is totally geodesic for this metric 
(conservation of linear momentum) and so we restrict the metric
to this linear subspace  $V \cong \C^2 \subset \C^3$.  

{\bf Step 2. Symmetry reduction.}  Observe that $ds^2_{JM}$ is invariant under complex scaling: $q \mapsto \lambda q,  \lambda \in \C^*: = \C \setminus \{0 \}$.
As a result, the metric  $ds^2 _{JM}$ descends to the quotient $\C^2 \setminus \{0 \} / \C^*$ which is $\C \bP ^1$, the shape sphere.
When we use this pushed-down metric the quotient projection  map $\V_{cm} \setminus \{0 \} \to \C \bP^1$ becomes a Riemannian submersion. 
We can identify this quotient projection with the shape projection $\pi$ composed with the radial projection
$w \mapsto \frac{1}{|w|} (w)$ from the shape space minus the origin onto a shape sphere.  
 The Riemannian submersion property implies  that geodesics for $ds^2 _{JM}$ which are orthogonal to the $\C^*$-fibers project
  project onto geodesics for the quotient metric.   Now use the facts that  $\dot I = 2  \langle q, v \rangle$, $J =  \langle v, i q \rangle$,
and   $ds^2 _{JM}$ is conformal to the mass metric $ds^2_E$, to conclude that a curve is   orthogonal to the fiber if and only if   $\dot I = 0 = J$ along that curve
to see that our Newton solutions are projected to shape space geodesics for the pushed-down metric.    It remains to compute  this pushed down metric,
let us call it $\pi_* ds^2 _{JM}$, and show it
is the Jemisphere metric.  
In \cite{me_monthly}, eq (43) it is shown that $ds^2 _E = \frac{|dw|^2}{I}$ when pushed down to the  shape {\it space}.  We have already seen that $\Delta = \mu w_3$
so that $U = \frac{\gamma \mu^2  I}{w_3^2}$. 
It follows that $\pi_* ds^2_{JM}= \frac{\gamma \mu^2 |dw|^2}{w_3 ^2}$, which is the Jemisphere metric, as required.  (The scaling constant $\gamma \mu^2$ changes
the constant curvature to $-1/\mu \sqrt{\gamma}$ but does not change the geodesics.) 

\begin{figure}[h]
\scalebox{0.4}{\includegraphics{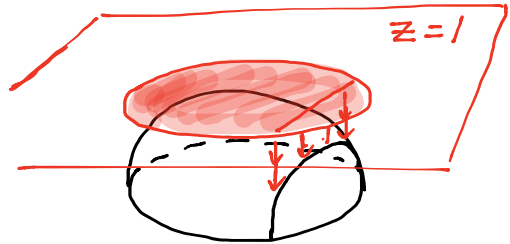}}
\caption{Projection along the vertical is an isometry between the Klein and Jemisphere models
of the Hyperbolic Plane.}  \label{Klein_to_Jemi}
\end{figure}

Finally we want to describe the geodesics for the Jemisphere metric in terms of our Hopf-Jacobi coordinates.
Use  the isometry between the Klein model and the Jemisphere metric which
is projection along   the $w_3$-axis, as described in \cite{All_Models}.  See figure \ref{Klein_to_Jemi}.
More specifically , the Klein model is realized as the disc $w_1 ^2 + w_2 ^2 < 1$ lying on the
plane $w_3 =1$ tangent to the sphere at the north pole.  (We've assumed $I_0 =4$ for simplicity so that $|w| =1$.) Project  points $(w_1, w_2, 1)$ of the Klein model
to  points $(w_1, w_2, w_3)$,  of the hemisphere $w_1 ^2 + w_2 ^2 + w_3 ^2 = 1$, $w_3 >0$ to get the isometry between models.  (Note $w_3 = \sqrt{ 1 -w_1 ^2 - w_2 ^2}$.)
Since   geodesics of the Klein model are chords $A w_1 + B w_2 = C$ and since  they map over to geodesics of the Jemisphere,
we have that the    same linear equations in $w_1, w_2$
characterize geodesics  in the Jemisphere model.

QED. 

\section{Confession. Motivation.  History. Open Problems} 

 I teach a  geometry class nearly every year for an audience consisting  primarily of  future high school mathematics teachers.    I feel an obligation to teach the rudiments of   hyperbolic geometry.
 But almost  all   the students leave the class  with no  understanding of what the hyperbolic plane is  beyond a vague sense  that ``things get really tiny and squinched together when you approach the
   x-axis'' .  The present paper
began as an attempt to provide a  natural road in to hyperbolic geometry and its intuition for these students.  I believe I  have failed.  Nevertheless I
hope some readers find the   exercise was interesting.

The papers \cite{hyperbolic_pants} and \cite{connor} combined  to suggest the  approach taken here to building the  ``mechanical hyperbolic metric'' of theorem 1. 
 In \cite{hyperbolic_pants} I applied  the reductions of  the present  
paper to the ``strong-force'' $1/r^2$ potential $V_2= -\Sigma {m_a m_b} {r_{ab}^2} $ in the equal mass three-body case to obtain a metric on
the shape sphere minus its three binary collision points, i.e. on the topologist's pair-of-pants.  The main theorem there
is that this metric is complete with negative Gauss curvature everywhere except at two points (the equilateral triangles of Lagrange).  As a corollary, the figure eight solution for that potential
is unique up to isometry and scaling.    In \cite{connor} , Connor Jackman and I  tried to  extend the $N=3$ hyperbolicity  to the
case $N=4$ body problem with  equal masses.  The potential   has the same form $V_2$ , except the sum is now taken  over all 6 pairs amongst the 4 bodies.
Connor proved that   the   results of \cite{hyperbolic_pants} do not readily extend: the resulting 4-body JM metric
has mixed curvature: some 2-planes have positive curvature, some negative. 

The tricks we applied in  \cite{hyperbolic_pants} ,  \cite{connor} and the present paper apply to any 
 Galilean-invariant potential of homogeneity $-2$ on the N-body configuration space $\C^N$, provided that potential $V$  is negative (possibly $-\infty$).
The result is a metric on $\C \bP ^{N-2} \setminus \Sigma$ where $\Sigma$ is the set on which $V = -\infty$.  That JM metric
has the form $U ds^2 _{FS}$  where $U = -V$.  
 
  Reflecting on the jump from N = 3 to N=4 in  jumping from \cite{hyperbolic_pants} to  \cite{connor}, it is natural to  wonder
what  might happen when we take such a   jump   for our  ``3-point''   potential (\ref{potential}).   To define the analogous N-body potential we 
sum our 3-point  potential over all  triples of bodies. 
Thus,  consider   N point masses  in the plane, labeling them $\{1, 2, \ldots, N \}$.   
For each choice of 3 indices $i,j,k$ out of $\{1, 2, \ldots, N\}$, let $I(i,j,k)$ be the moment of inertia of the triangle
 formed by vertices $q_i, q_j, q_k$  and   let $\Delta (i,j,k)$ be the signed area of this triangle.  Then the N-body analogue of our potential is 
$$V_N  = - \gamma \Sigma_{i,j,k} \frac{I(i,j,k)}{\Delta (i,j,k)^2}$$
the sum being over all three-element subsets  $\{i,j,k \} \subset \{1, 2, \ldots, N\}$ .  We observe that the singularity of $V_N$ is precisely the set $\Sigma_N$ of all  ``non-generic' planar N-gons,
where we say that an N-gon is non-generic if some 3 of its vertices are collinear.   $\Sigma_N$ is a real codimension 1 subvariety
which cuts the configuration space into a number of components.  

{\bf Question 1.}  For $N=4$ is the resulting JM  metric on $\C \bP^2 \setminus \Sigma_4$ hyperbolic?

\vskip .3cm 
One computes without much difficulty that this  JM metric is complete. Indeed, near
any typical point of $\Sigma_N$ we can choose a coordinate $w_3$ such that $w_3 = 0$ locally
defines $\Sigma$ and the JM metric as $w_3 \to 0$ has  leading asymptotics  that of the
hyperbolic metric:  $\frac{ 1}{w_3 ^2} ds^2 _{Euc}$.   There are restrictions on the topology of complete hyperbolic manifolds. 
If the manifold is simply connected then it must be diffeomorphic to the ball in that dimension.     This suggests looking into the
topology of $\C \bP^2 \setminus \Sigma_4$.    We have verified that $\C \bP^2 \setminus \Sigma_4$ consists of 14 components, and that 
each one of which is diffeomorphic to $B^4$, providing weak circumstantial evidence that the answer might be ``yes'' to question 1.

Jumping further ahead to higher $N$:

{\bf Question 2. }  How many components are there in the space $\C \bP^{N-2} \setminus \Sigma_N$ of general position planar  N-gons?
Is each component diffeomorphic to an open ball of dimension $2(N-2)$?

\vskip .3cm

{\bf Question 3.} Do potentials of the form of eq (\ref{potential}) arise in any physical or chemical problems of 
interest?

\section{Acknowledgements}

I am grateful to Persi Diaconis and Toshiaki Fujiwara for inspirational conversations, and to 
Serge Tabachnikov and Richard Schwartz for conversations arount the open questions.   
 I  thank NSF grant DMS-1305844 for essential support.

\end{document}